\documentclass[10pt]{article}
\usepackage{amssymb,amsmath,latexsym,epsfig,amsfonts}
\usepackage{amsthm,amscd}
\newcommand{\mgn}{\overline{\mathcal{M}}_g^n}
\makeatother
\newtheorem{theorem}{Theorem}[section]
\newtheorem{corollary}[theorem]{Corollary}
\newtheorem{lemma}[theorem]{Lemma}
\newtheorem{proposition}[theorem]{Proposition}
\newtheorem{remark}[theorem]{Remark}
\newtheorem{example}[theorem]{Example}

\newtheorem{definition}[theorem]{Definition}

\def\proof{\noindent{{\bf Proof:}}\,\,} \def\qed{{\hfill{\large
$\Box$}}}

\title{Euler Characteristics of Moduli Spaces of Curves
\thanks{Research by the first author was partially supported by MIUR and GNSAGA. Research by the second author was partially supported by
NSF under grant DMS-01-07621.}}

\begin{document}

\author{Gilberto Bini\thanks{Dipartimento di Matematica
            ``F. Enriques'', Universit\`{a} degli Studi di Milano,
            Milano, Italy. E-mail: {\tt bini@mat.unimi.it}}, John
            Harer\thanks{Department of Mathematics, Duke University,
            Durham, North Carolina, (U.S.A.). E-mail: {\tt
            John.Harer@Duke.edu}}}

\maketitle

\begin{abstract}
Let ${\mathcal M}_g^n$ be the moduli space of $n$-pointed Riemann
surfaces of genus $g$. Denote by $\mgn$ the Deligne-Mumford
compactification of ${\mathcal M}_g^n$. In the present paper, we
calculate the orbifold and the ordinary Euler characteristic of $\mgn$
for any $g$ and $n$ such that $n > 2-2g$. 

\end{abstract}

\section{Introduction}\label{intro}
The moduli space of $n$-pointed Riemann surfaces of genus $g$,
${\mathcal M}_g^n$, is an object of much importance in several branches
of mathematics and theoretical physics.  It parameterizes algebraic
curves of genus $g$ with $n$ points, or equivalently Riemann surfaces
of genus $g$ with $n$ points.  Mumford and Deligne in \cite{Mum}
defined a natural compactification of ${\mathcal M}_g^n$, $\mgn$, by
adjoining stable curves at infinity.  These spaces serve as
classifying spaces in algebraic geometry, so it is very important to
understand their topological structure, especially their homology and
cohomology.  Madsen \cite{Ma} proved the Mumford conjecture by calculating the stable cohomology of ${\mathcal M}_g^n$.  In the current paper we work in the unstable range, calculating both the
orbifold and the ordinary Euler characteristics of $\mgn$ for any $g$
and $n$ with $2g - 2 + n > 0$.

\paragraph{Related work.}
In \cite{Haz}, Harer and Zagier calculated the orbifold Euler
characteristic of ${\mathcal M}_g^n$. Moreover, they computed the
ordinary Euler characteristic of ${\mathcal M}_g^n$ for any $g$ and
$n=0,1$. Since then, there have been some results on the ordinary
Euler characteristic of $\mgn$ only for low values of $g$. The reader
is referred to \cite{Ge} and \cite{YM} for $g=0$, to \cite{Ge1} for
$g=1$, to \cite{bgm}, \cite{Fa}, \cite{Ge2}, for $g=2$ and to
\cite{Ge3} for $g=3$.

\paragraph{Results.}
The main results of this paper are formulae for the Euler
characteristics of ${\mathcal M}_g^n$ and $\mgn$.  These are given in
theorems $\ref{euler1}$, $\ref{statest}$ and $\ref{matform}$.  These
theorems are proven by combining the techniques of generating
functions, integral representations and Wick's lemma that are well
established in the field with other counting methods, including a
formula of Serre and Brown \cite{Bro} (also used in \cite{Haz}) for
computing actual Euler characteristics in terms of orbifold ones.  We
compute tables of values in all cases (these were generated using
Maple) and present the values below.  Our results agree with previous
computations made in the papers mentioned above.

\paragraph{Outline.} 
Section \ref{sec1} describes the stratification of $\mgn$ in terms of
stable curves.  In \ref{sec2} we review the results of \cite{Haz}, and
compute the Orbifold Euler characteristic of ${\mathcal M}_g^n$.
Section \ref{sec4} concludes the paper by calculating the actual Euler
characteristic of ${\mathcal M}_g^n$. Tables of values are given in
each case.

\section{The Graph-type Stratification of ${\overline{\cal{M}}}_g^n$}
\label{sec1}

In this section, we review some basic facts and
definitions we shall use in the rest of the paper.  The
moduli space of stable curves ${\overline{\cal{M}}}_g^n$
admits a stratification, which
is determined by the configuration of nodes and
by irreducible components of $n$-pointed genus $g$
stable curves. This stratification may be described via
stable graphs. For the sake of completeness, we briefly
recall their definition. 

\begin{definition}
\label{stablegr}
Let $g,n$ be non-negative integers such that $n > 2-2g$.
A stable graph of type $(g,n)$ is given by the following
data:

\begin{itemize}
\item [{\bf SG1)}] two finite sets $V_G$ and $L_G$;
\item [{\bf SG2)}] a partition ${\mathcal P}$ of $L_G$ into subsets with one or two elements;
\item [{\bf SG3)}] a map $\gamma$ from $V_G$ to the set of integers $\{0,\ldots,g\}$ such that
$$
g=\sum_{v \in V_G}\gamma(v) +h^1(G);
$$
\item [{\bf SG4)}] a subset $L(v) \subset L_G$, $ v \in V_G$, such that $2\gamma(v)-2+|L(v)| >0$;
\item [{\bf SG5)}] a map $\nu$ from subsets of $L_G$ with one element
to $\{1, \ldots, n\}$.
\end{itemize}

\end{definition}

The elements of $V_G$ are the vertices of $G$, whereas
the elements of $L_G$ are the half-edges of $G$.
Moreover, we shall call {\em legs} the subsets of
$\mathcal{P}$ with one element and {\em edges}
the subsets of $\mathcal{P}$ with two elements. {\bf SG3)}
relates $g$ to the structure of $G$. In fact,
$h^1(G)=1-v(G)+e(G)$, where  $v(G)$ and $e(G)$
denote the number of vertices and edges of a stable graph,
respectively. The automorphism group $Aut(G)$ of $G$ is the set of
bijections which map $V_G$ to $V_G$, $L_G$ to $L_G$, and
preserve all data of Definition \ref{stablegr}.
In what follows, for the sake of simplicity,
we shall denote $V_G$ and $L_G$ by $V$ and
$L$, respectively. Moreover, by abuse of notation,
we will call $\gamma(v)$ the genus of $v$, and $g$
the genus of $G$.

Given a stable graph $G$ of type $(g,n)$, choose
an ordering of $L(v)$ for each vertex $v$. Next, consider
the morphism \begin{equation}
\label{csi}
\xi _{G }:\prod _{v\in V}\overline{\mathcal {M}}_{\gamma(v)}^{l(v)}\rightarrow
\overline{\mathcal {M}}_g^n,
\end{equation}
where $l(v)=|L(v)|$. A point in the domain is the datum of an
$l(v)$-pointed curve $C_v$ for each $v$.  The image point is the
$n$-pointed genus \( g \) curve which is obtained as follows: identify
the marked points of $C_v$ corresponding to the half-edges of $G$
which are connected by an edge. By definition, the map $\xi _{G }$ is
independent from the ordering of the sets $L(v)$'s. Set, further,
$\Delta^o_G=\xi^o({\cal{M}}_G)$ and denote by $\Delta_G$ its closure
in $\overline{\mathcal {M}}_g^n$. By definition of (\ref{csi}), two
elements in a fiber of $\xi_G^o$ differ by an automorphism of $G$.
This means that

\begin{equation}
\label{edamojelo}
\Delta_G^o \simeq {\cal{M}}_{G}/Aut(G).
\end{equation}

We recall that a stable graph $G$ of type $(g,n)$
degenerates to a stable graph $G'$ of the same type
if $G'$ can be obtained from $G$ by a chain of the
following moves:
\begin{enumerate}
\item collapse an edge that joins two different
vertices $v_1$ and $v_2$ and label the new vertex $v$ with
$\gamma(v)=\gamma(v_1)+\gamma(v_2)$,
\item collapse a loop to a vertex $v$ and increase the
genus of $v$ by one.
\end{enumerate}

In any case, we write $G'<G$. Thus the following holds:

\begin{equation}
\label{strat}
\overline{\cal{M}}_g^n=\bigcup_G \Delta_G^o,
\end{equation}
where the union is over stable graphs of type $(g,n)$.

The locally closed strata $\Delta_G^o$ are equipped with an
orbifold structure. Fix a topological oriented surface
$S_{\gamma(v), l(v)}$ for each vertex $v$ of $G$. Such a surface
has genus $\gamma(v)$ and $l(v)$ marked points,
where $l(v)$ is the number of half-edges outgoing from $v$.
Denote by $a(v)$ the number of half-edges outgoing from $v$
that are legs, and set
$b(v):=l(v)-a(v)$. For each vertex, consider the
Teichm\"{u}ller space $\mathcal{T}_{\gamma(v)}^{l(v)}$
and the mapping class group
$\Gamma_{\gamma(v)}^{l(v)}$. If

\begin{equation}
\label{latidigi}
T(G):= \prod_{v \in V_G}\mathcal{T}_{\gamma(v)}^{l(v)},
\end{equation}
the orbifold structure of $\Delta_G^o$ can be described as
follows. We recall that the elements of $Aut(G)$ are obtained
as compositions of permutations of vertices $v_1$ and $v_2$ (when
$\gamma(v_1)=\gamma(v_2)$, $a(v_1)=a(v_2)=0$, $b(v_1)=b(v_2)$)
or permutations of half-edges of $G$. Any permutation of
the $l(v)$ half-edges outgoing from $v$ induces a permutation of
the $l(v)$ marked points of $S_{\gamma(v), l(v)}$.
Accordingly, $Aut(G)$ acts on $T(G)$ as follows.
Fix a vertex $v$ and, for simplicity, denote by
$\left[C; x_{1}, \ldots, x_{l}, \left[f \right]\right],
l=l(v)$, an element of the Teichm\"{u}ller space associated with
$v$. Any $\tau$ in $Aut(G)$, which permutes
the half-edges outgoing from $v$, maps
$\left[C;x_{1}, \ldots, x_{l},[f]\right]$ to
$\left[C;x_{\tau(1)}, \ldots, x_{\tau(l)},[\tau \circ f \circ
\tau^{-1}]\right]$ - $\tau$ is the permutation of the marked
points of $S_{\gamma(v),l(v)}$. On the other hand,
take $\tau$ in $Aut(G)$
such that $\tau(v_1)=v_2$, where $v_1$ and $v_2$ are vertices of $G$
with the same genus and $a(v_1)=a(v_2)=0$, $b(v_1)=b(v_2)$. Set
$l_1=l(v_1)=l(v_2)$, and consider the elements
\begin{equation}
\label{duelementi}
\left[C_1; x_{1}, \ldots, x_{l_1},[ f_1]\right], \, 
\left[C_2;y_{1}, \ldots, y_{l_1},[f_2]\right].
\end{equation}
These two are exchanged by $\tau$ in the product $T(G)$.

\begin{definition}
\label{arieccote}
\begin{equation}
\Gamma(G):= \prod_v \Gamma_{\gamma(v)}^{l(v)} \rtimes Aut(G)
\end{equation}
\end{definition}

We now define an action of $Aut(G)$ on $\prod_v
\Gamma_{\gamma(v),l(v)}$ so that $\Gamma(G)$ acts
on $T(G)$ and $T(G)/\Gamma(G)\cong \Delta_G^o$.
If $\tau \in Aut(G)$ permutes the half-edges of $G$ and
$ \Bigl(\prod_v \left[
{\mathfrak h}_v\right] \Bigr)\in \prod_v \Gamma_{g(v)}^{l(v)}$, set
\begin{equation}
\label{conimezzilati2}
\tau \cdot \Bigl(\prod_v \left[ {\mathfrak h}_v\right] \Bigr)= \Bigl(\prod_v \left[\tau \circ {\mathfrak h}_v \circ \tau^{-1}\right] \Bigr).
\end{equation}

On the other hand, if $\tau \in Aut(G)$ permutes two vertices $v_1$
and $v_2$ with the same genus and $a(v_1)=a(v_2)=0$, $b(v_1)=b(v_2)$,
consider two elements as in (\ref{duelementi}). If ${\mathfrak h}_1$
and ${\mathfrak h}_2$ are elements in the mapping
class groups associated with $v_1$ and $v_2$, then
$\tau$ acts on the group $\Gamma(G)$ since it swaps
${\mathfrak h}_1$ and ${\mathfrak h}_2$. As a result, the
semi-direct product $\Gamma(G)$ in (\ref{arieccote}) is
well defined and acts on $T(G)$ in the following way.
First, let us consider the case of
an automorphism $\tau$ which permutes the half-edges of $G$
that stem from a vertex $v$. If $\left[C;x_1,
\ldots, x_{l}, \left[f\right]\right]$ and ${\mathfrak h}$
belong to the Teichm\"{u}ller space and to
the mapping class group associated
with $v$, $\left[C;x_{1}, \ldots,x_{l},[f]\right]$
is mapped to $\left[C;x_{\tau(1)}, \ldots,x_{\tau(l)},[\tau \circ
{\mathfrak h} \circ f \circ \tau^{-1}]\right]$. Second,
let $\tau$ permute two vertices of $G$, $v_1$ and $v_2$,
with no legs, the same genus, and the same number of half-edges.
Consider two elements as the ones in
(\ref{duelementi}), which belong to the Teichm\"{u}ller spaces
corresponding to $v_1$ and $v_2$, and two classes ${\mathfrak h}_1$
and ${\mathfrak h}_2$ in the mapping class groups associated with
$v_1$ and $v_2$. Then the action of $\Gamma(G)$ exchanges $
\left[C_1;x_1, \ldots, x_{l_1}, \left[{\mathfrak h}_1\circ
f_1\right]\right] $ with $ \left[C_2;y_1, \ldots, y_{l_1},
\left[{\mathfrak h}_2\circ f_2\right]\right]$.

It is easy to check that the elements in $T(G)/\Gamma(G)$
are obtained by looking at the orbit of pointed stable
curves under the action of $Aut(G)$. Thus,
\begin{equation}
\label{quot}
T(G)/\Gamma(G)\cong \Delta_G^o.
\end{equation}

Furthermore, $\Delta_G^o$ has an orbifold structure, since the action
of $\Gamma(G)$ is properly discontinuous and with finite stabilizers,
as can be readily checked.

\section{The Orbifold Euler Characteristic of $\mgn$}\label{sec2}

The orbifold structure of $\mgn$ naturally induces
an orbifold Euler characteristic, which will be hereafter
denoted by $\chi(\mgn)$. In this section, we use the stratification
described in (\ref{strat}) to determine generating functions of
the rational numbers $\chi(\mgn)$.

Suppose an orbifold $M$ admits a manifold $\widetilde{M}$
as a finite branched covering $\pi:{\widetilde M}
\rightarrow M$ of degree $d$. Then $\chi(M)$ turns out to
be $e(\widetilde M)/d$, where $e(\widetilde M)$ is
the ordinary Euler characteristic of ${\widetilde M}$.
Recall that the Euler characteristic of a virtually torsion free group $H$ is defined
similarly, i.e., $ \chi(H)=\chi(\widetilde{H})/d$, where
$\widetilde{H}$ is a torsion free subgroup of index $d$ in
$H$. We shall use this group theoretic analogy to
compute $\chi(\Delta^o_G)$.

First, observe that $\Gamma(G)$ contains torsion free
subgroups $\widehat{\Gamma}(G)$ which act freely on $T(G)$.
This follows from well known facts about level structures of algebraic curves. 
As a consequence, $T(G)/\widehat{\Gamma}(G)$ is a
finite branched covering of $\Delta^o_G$ of degree
$[\Gamma(G):\widehat{\Gamma}(G)]$. Therefore, $
\chi(\Delta^o_G)=\chi(\Gamma(G))$. By the short exact sequence of
groups
$$
1 \rightarrow \prod_v \Gamma_{\gamma(v)}^{l(v)} \rightarrow \Gamma(G) \rightarrow Aut(G) \rightarrow 1,
$$
we get 
\begin{equation}
\label{lascomponi}
\chi(\Gamma(G))=\left(\prod_v\chi({\cal
M}_{\gamma(v)}^{l(v)})\right)/|Aut(G)|.
\end{equation} Thus, $$
\chi(\overline{\mathcal M}_g^n) = \sum_G
\frac{\prod_v\chi\left({\cal M}_{\gamma(v)}^{l(v)}\right)}{|Aut(G)|}.$$

The orbifold Euler characteristic of the moduli space
${\mathcal M}_g^n$  has been computed in \cite{Haz}.
More precisely, the following holds.

\begin{theorem} (\cite{Haz}) For non-negative integers $g$, $n$, $n >2-2g$, the orbifold Euler characteristic of ${\mathcal M}_g^n$ is
$$
\chi({\cal{M}}_g^n)=(-1)^n\frac{(2g-1)B_{2g}}{(2g)!}(2g+n-3)!,
$$
where $B_{2g}$ is the $(2g)$-th Bernoulli number.
\end{theorem}

In order to compute $\chi({\overline{\mathcal M}}_g^n)$, we
introduce the power series
$$
F(x, \hbar):= \sum_{g \geq 0}F_g(x)\hbar^{g-1},
$$
where 
$$
F_g(x):=\sum_{n > 2-2g, n \geq 0} \chi(\overline{\cal{M}}_g^n)\frac{x^n}{n!}.
$$

We will express the formal power series $F(x, \hbar)$ in terms of the known generating series
$$
\Omega(x,\hbar):= \sum_{g \geq 0} \sum_{\underset{n \geq 0}{n > 2-2g}} \chi({\cal{M}}_g^n)\frac{x^n}{n!}.
$$
Standard techniques in asymptotic theory will
yield closed formulas for $F_g(x)$.

\begin{theorem}
\label{euler1}
\begin{equation}
\label{wickone}
\exp(F(x, \hbar))=\int_{\mathbb{R}}\exp{\Bigl(-\frac{(x-y)^2}{2\hbar}+\Omega(y,\hbar)\Bigr)}\frac{dy}{\sqrt{2\pi \hbar}}.
\end{equation}
\end{theorem}
{\bf Proof.} If we make the substitution $y-x=z\sqrt{\hbar}$, the integral on the right-hand side of (\ref{wickone}) reduces to a one-dimensional gaussian integral which can be computed directly. Moreover, if the exponential to be integrated is expanded as a power series, we get
\begin{eqnarray}
\label{svi1}
1 &+&   \sum_{k\geq 1}\sum_{g_1,\ldots,g_k \geq 0}
        \sum_{\underset{ r_j >2-2g_j}{r_1, \ldots, r_k}}
        \prod_{j=1}^k\chi({\cal{M}}_{g_j,r_j}) \cdot \\
&\cdot& \sum_{\underset{\sum t_j even}{t_1, \ldots t_k=0}}^{r_1, \ldots, r_k}
        \frac{(t_1 + \ldots + t_k -1)!!}{k!t_1! \ldots t_k!}
        \frac{x^{\sum_{j=1}^k(r_j-t_j)}}{\prod_j (r_j-t_j)!}
        \hbar^{\sum_{j=1}^k(g_i-1)+ \frac12 \sum_j t_j}. \nonumber
\end{eqnarray}

The claim will follow if the sum in (\ref{svi1}) can be
rewritten as a sum over stable graphs. For this purpose,
consider $k$, $ k \geq 1$, stable graphs $G_1,\ldots, G_k$ each of which has one vertex of genus
$g_j$ and $r_j$ legs. If we choose $t_j$ legs, $0 \leq t_j \leq r_j$,
from each $G_j$, there are $(t_1+\ldots +t_k-1)!!$ possible ways of
interconnecting them, provided $\sum t_j$ is even. Such
pairing yields a disconnected stable graph $G$ of type $(g_k,n_k)$, where
$$
g_k=\sum_{j=1}^kg_j+1-k+\frac12\sum_{j=1}^kt_j, \quad n_k=\sum_{j=1}^k(r_j-t_j).
$$

Contrarily, if we fix non-negative integers $g$ and $n$ ($n > 2-2g$) and a disconnected stable graph of type $(g,n)$, we can determine a collection of integers $k, t_1, \ldots, t_k, r_1, \ldots, r_k$ as in the sum which appears in (\ref{svi1}). This sum can therefore be rearranged as 
\begin{equation}
1+\sum_{g \geq 0}
\sum_{\substack{n \geq 0 \\ n > 2-2g}}
\sum_{G \in \mathcal{G}_{g,n}}
\frac{\chi({\cal{M}}_G)}{|Aut(G)|}\hbar^{g-1},
\label{svi2}
\end{equation}
where $\mathcal{G}_{g,n}$ is the set of disconnected stable graphs of
genus $g$ with $n$ legs. By standard combinatorial arguments, the
theorem is completely proved.
\begin{flushright}
$\Box$
\end{flushright}

\subsection{Asymptotic Formulas For $F_g(x)$}\label{sec21}

In order to deduce formulas for $F_g(x)$ we perform a {\em
semiclassical expansion} of the integral on the right-hand side of
(\ref{svi1}). In other
words, we substitute
$$
U(x,y, \hbar):= -\frac{(x-y)^2}{2\hbar}+\Omega(y,\hbar)
$$
with its formal power series centered at the solution
of
\begin{equation}
\label{sap}
\overline{y}=x+ \sum_{g \geq 0} \Omega_g'(\overline{y})\hbar^g,
\end{equation}
where the prime denotes the derivation with respect to
the variable $y$. We thus look for a solution of 
(\ref{sap}) of the form
$$
\overline{y}(x, \hbar):=\sum_{g \geq 0}y_g(x)\hbar^g.
$$
This yields the recursive relations
\begin{equation}
\label{momeserve}
y_0(x)=x+\sum_{n \geq 2}\chi({\cal{M}}_{0,n+1})\frac{y_0^n(x)}{n!},
\end{equation}

\begin{equation*}
y_g(x)= \sum_{s=0}^g\sum_{\underset{n \geq 0}{n > 1-2s}}
\chi({\cal{M}}_{s,n+1})\sum_{\substack {m_1+2m_2+ \ldots +gm_g=g-s \\
                              m_0+ m_1+\ldots +m_g=n}}
\frac{y_0^{m_0}(x) \ldots y_g^{m_g}(x)}{m_0! \ldots m_g!}.
\end{equation*}

The function $y_0(x)$ can be computed via the differential equation
\begin{equation*}
\frac{dy_0(x)}{dx}(1-\log(1+y_0(x)))=1.
\end{equation*}
This yields the power series
\begin{equation}
\label{zapatos}
y_0(x)=x+\frac{x^2}{2}+\frac{x^3}{3}+\frac{7}{24}x^4+\frac{17}{60}x^5 +
\frac{71}{240}x^6 +\frac{163}{504}x^7 +o(x^8).
\end{equation}

Since $y_0'(0)=0$, the $y_g(x)$'s are uniquely defined via the recursive relations
\begin{equation*}
\frac{y_g(x)}{y_0'(x)}=\sum_{s=1}^g\sum_{\substack{n >1-2s \\ n \geq 0}}
\chi({\cal{M}}_{s,n+1})\cdot \sum_{\substack{m_1 +2m_2+ \ldots + gm_g=g-s \\ m_0+m_1+ \ldots +m_g=g}}
\frac{y_0^{m_0}(x) \ldots y_g^{m_g}(x)}{m_0! \ldots m_g!},
\end{equation*}

Let us now expand the function $U(x,y,\hbar)$ about the point
$\overline{y}(x, \hbar)$, and set $w=y-\overline{y}$. Thus, we get
$$
-\frac{(x-\overline{y})^2}{2\hbar}+\Omega(\overline{y},\hbar)
-\frac{1}{2\hbar}w^2\bigl(1-\sum_{g \geq
0}\Omega^{(2)}_g(\overline{y})\hbar^{g}\bigr) +\sum_{k \geq
3}\frac{1}{k!}w^k\Bigl(\sum_{g \geq 0}
\Omega_g^{(k)}(\overline{y})\hbar^{g-1}\Bigr),
$$
where the superscript $(j), j \geq 2$, denotes the derivation with
respect to the variable $y$. For the sake of simplicity, we set
$$
G(x, \overline{y}(x, \hbar)) = \sum_{g \geq 0} \Omega_g ^{(2)}(\overline{y}(x,\hbar))\hbar^g,
$$
$$
S_{k}(x, \overline{y}(x, \hbar))= \sum_{g \geq 0} \Omega_{g}^{(k)}(\overline{y}(x, \hbar)) \hbar^{g-1},
$$
and
\begin{eqnarray*}
A(x, \overline{y}( x, \hbar))& = &\sum_{r \geq 1} \sum_{\underset{\sum k_i \, even}{k_1, \ldots , k_r \geq 3}}
\frac{(k_1 + \ldots k_r-1)!!}{k_1! \ldots k_r!} \cdot \\
&\cdot&\frac{S_{k_1}(x,\overline{y}(x,\hbar))\cdot  \ldots \cdot S_{k_r}(x,\overline{y}(x,\hbar))}
{\sqrt{1-G(x,\overline{y}(x,\hbar))^{k_1 + \ldots k_r+1}}} \hbar^{\frac12 \sum{k_i}}.
\end{eqnarray*}

Then the following holds.
\begin{theorem}
\label{asfi}
Let
$$
F(x, \hbar):=\sum_{g \geq 0} \sum_{\substack{n \geq 0 \\ n> 2-2g}}
\chi(\overline{\cal{M}}_g^n)\frac{x^n}{n!}\hbar^{g-1}
$$
be the generating function of $\chi(\overline{\cal M}_g^n)$.  An
asymptotic expansion of $F(x,\hbar)$ is given by
\begin{eqnarray}
& &\frac{(-x + \overline{y}(x,\hbar))^2}{2 \hbar} + \sum_{g \geq 0} \Omega_g(\overline{y}(x,\hbar))\hbar^{g-1} \label{asfiu} \\
&-&\frac12 \log(1-G(x,\overline{y}(x,\hbar)) \label{asfid} \\
&+& \log(1 +A(x, \overline{y}(x, \hbar)). \label{asfit}
\end{eqnarray}
\end{theorem}

{\bf Proof.} The claim follows by rewriting the integral appearing on the right-hand side of (\ref{svi1}) as
$$
\exp{\Bigl(\frac{(-x + \overline{y}(x,\hbar))^2}{2 \hbar} + \sum_{g
\geq 0} \Omega_g(\overline{y}(x,\hbar))\hbar^{g-1}\Bigr)} \cdot
$$
\begin{equation}
\cdot \int_{\mathbb{R}} \exp
{\left(-\frac{1}{2\hbar}w^2\Bigl(1-\sum_{g \geq
0}\Omega^{(2)}_g(\overline{y}) \hbar^{g}\Bigr)\right)}
\exp{\left(\sum_{k \geq 3}\frac{1}{k!}w^k\Bigl(\sum_{g \geq 0} \Omega_
g^{(k)}(\overline{y})\hbar^{g-1}\Bigr)\right)}\frac{dw}{\sqrt{2\pi
\hbar}} \label{esto}
\end{equation}

By Wick's Lemma (see \cite{Biz}), note that the term involving $A(x,
\hbar)$ originates from the expansion of the exponential in the
integral (\ref{esto}).
\begin{flushright}
$\Box$
\end{flushright}

The semiclassical expansion used in Theorem \ref{asfi} can be
interpreted as a {\em loop-wise expansion}, i.e., as an expansion with
respect to the first Betti number of a stable graph of type
$(g,n)$. Thus, we shall describe the
$y_g(x)$'s, $g \geq 0$, from a combinatorial point of view.

\begin{definition}
\label{dajeunnome}
A stable tree $T$ of genus $h, h \geq 0,$ is a tree such that
\begin{itemize}
\item [T1] there exists a map $\gamma : V_T \rightarrow \{0, \ldots, h\}$ with $\sum_{v}\gamma(v)=h$,
\item [T2] there are $n$, $n \geq 0$, numbered leaves and $j$ unnumbered leaves going into the root such that $j \geq 1$, $n +j>2 -2h$,
\item [T3] for every vertex $v$ of the tree,
the number of outgoing edges (including the leaves) is greater than
$2-2\gamma(v)$.
\end{itemize}
\end{definition}

In the sequel, we shall denote by $\mathcal{G}'_{h, n+j}$ the
collection of stable trees of genus $h$, with $n$ numbered
leaves and $j$ unnumbered leaves going into the root.  A graph in $\mathcal{G}_{h,n+j}'$ is by definition a stable graph of type $(h, n+j)$. Take now the two generating functions
$$
\xi_0(x):=x+\sum_{n \geq 2}\sum_{T \in {\mathcal{G}}_{0,n+1}'}\chi(\Delta_T^o)\frac{x^n}{n!},
$$
$$
\xi_{g}(x):=\sum_{h=0}^g\sum_{n > 1-2h}\sum_{T \in {\mathcal{G}}'_{h,n+1}}
\chi(\Delta_T^o)\frac{x^n}{n!},
$$
where $\chi(\Delta_T^o)$ is the orbifold Euler characteristic of the
open stratum defined as the image of the morphism $\xi_{T}^o$ in
(\ref{edamojelo}).

\begin{proposition}
\label{idc}
For each $g \geq 0$, $\xi_g(x)=y_g(x)$.
\end{proposition}

{\bf Proof.} Choose a graph $T \in \mathcal{G}'_{h,n+1}$. The root of $T$
corresponds to a vertex $v$ with an unnumbered leg, $t$ outgoing edges,
and $\gamma(v)=h, 0 \leq  h \leq g$. If we cut $T$ along these edges, we
get $m$ graphs $T_1, \ldots, T_m$ each of which belongs to
$\mathcal{G}'_{a,b+1}$, with $a \in \{0, \ldots, g-h\}, b>1-2a$.
The claim follows by (\ref{lascomponi}).
\begin{flushright}
$\Box$
\end{flushright}

Let ${\mathbb G}_{(g,l,n)}$ be the set of stable graphs with genus $g$,
$h^1(G)=l$, and $n$ legs. Theorem \ref{asfi} can be interpreted in the following way.

\begin{proposition}
\label{new}
\begin{itemize}

\item [i)] 
$$\frac{(-x + \overline{y}(x,\hbar))^2}{2 \hbar} + \sum_{g \geq 0}
\Omega_g(\overline{y}(x,\hbar))\hbar^{g-1} =\sum_{g \geq 0} \sum_
{\substack{n > 2-2g \\ n \geq 0}} \, \sum_{G \in {\mathbb G}_{(g,0,n)}}
\chi(\Delta^o_G) \frac{x^n}{n!}\hbar^{g-1};
$$
\item [ii)] $$ 
-\frac{1}{2}\log(1-G(\overline{y},\hbar))= \sum_{g \geq
1} \sum_ {\substack{n > 2-2g \\ n \geq 0}}\,\, \sum_{G \in {\mathbb
G}_{(g,1,n)}} \chi(\Delta^o_G) \frac{x^n}{n!}\hbar^{g-1};
$$

\item [iii)] 
$$
\log(1+A(x,\overline{y}(x,\hbar))) = \sum_{l \geq 2} \, \sum_{g \geq
l}\, \, \sum_{\substack{n > 2-2g \\ n \geq 0}} \, \, \sum_{G \in
\mathbb{G}_{(g,l,n)}} \chi(\Delta_G^o)\frac{x^n}{n!}\hbar^{g-1}.
$$
\end{itemize}
\end{proposition}

{\bf Proof}. i) Since 
$$
-\frac{(x-\overline{y})^2}{2\hbar}=-\frac{1}{2\hbar}\Bigl(\sum_{g \geq 0}\Omega_{g}'(\overline{y})\hbar^{g-1}\Bigr)^2,
$$
each contribution in (\ref{asfiu}) is of the form
$$
\chi({\cal{M}}_{h,r})P(x),
$$
where $P(x)$ is a polynomial in $y_{g}(x), g \geq 0$. By Proposition
\ref{idc}, this is just a sum of the Euler characteristics
$\chi(\Delta_G^o)$, where $G$ is a stable graph and $h^1(G)=0$.

ii) Observe that 
\begin{equation}
\label{derj}
\Omega_g^{(j)}(\overline{y})=\sum_{\substack{n \geq 0 \\ n > 2-j -2g}}\chi({\cal{M}}_{g,n+j})\frac{\overline{y}^n}{n!}.
\end{equation}

By using (\ref{sap}), we can describe (\ref{derj}) as
a sum over oriented rooted trees of arbitrary genus
having a root with $j$ unnumbered leaves. Therefore,
the product
$$
\Omega_h^{(j)}(\overline{y})\Omega_t^{(k)}(\overline{y})
$$
is a sum over stable graphs, which is obtained
in the following way. We
match the $j$ unnumbered leaves outgoing from the root
of $T_1 \in \mathcal{G}'_{a,j+n_1}, a \leq h,$ with the $k$
unnumbered leaves from the root of
$T_2 \in \mathcal{G}'_{b,k+n_2},b \leq t, n_1+n_2=n$.
In other words, multiplying derivatives of
$\Omega_g(\overline{y})$ gives rise to a sum over stable graphs with $h^1(G)
\geq 1$. In particular, since in (\ref{asfid}) there are only second
derivatives, the contribution
$$
-\frac{1}{2}\log(1-G(\overline{y},\hbar))
$$
can be rewritten as
$$
\sum_{g \geq 1} \, \, \sum_{\substack{n > 2-2g \\ n \geq 0}} \, \,
\sum_{G \in \mathbb{G}_{(g,1,n)}}\chi(\Delta_G^o)\frac{x^n}{n!}\hbar^{g-1}.
$$

iii) Analogously to the case $h^1(G)=1$, the contribution
in (\ref{asfit}) can be interpreted as
$$
\sum_{l \geq 2} \, \sum_{g \geq l}\, \, \sum_{\substack{n > 2-2g \\ n \geq 0}}
\, \,  \sum_{G \in \mathbb{G}_{(g,l,n)}}
\chi(\Delta_G^o)\frac{x^n}{n!}\hbar^{g-1}.
$$
\begin{flushright}
$\Box$
\end{flushright}

By Proposition \ref{new} we also have

\begin{theorem}
\label{domanisivede}
Let ${\cal{M}}_{g,n}^c$ be the moduli space of stable genus $g$ curves
with $n$ marked points whose associated stable graph $G$ is a
tree. Then
$$
\sum_{g \geq 0}\sum_{\substack{n \geq 0 \\ n >2-2g}}
\chi({\cal{M}}^c_{g,n})\frac{x^n}{n!}\hbar^{g-1}
=\frac{(-x + \overline{y}(x,\hbar))^2}{2 \hbar} + \sum_{g \geq 0}
\Omega_g(\overline{y}(x,\hbar))\hbar^{g-1}.
$$
\end{theorem}

\begin{example} {\rm The combinatorial interpretation carried out in Proposition \ref{new} yields explicit formulae for the functions $F_g(x)$. In the genus zero case,
\begin{equation*}
F_0(x)=\Omega_0(y_0)-\frac12\Omega'_0(y_0),
\end{equation*}
where $y_0(x)$ is defined in (\ref{zapatos}). Since the function $y_0(x)$ satisfies the identity
$$
(1+y_0)\log(1+y_0)=2y_0-x,
$$
the power series $x+F'_0(x)$ coincides with the one
given in \cite{YM} in an implicit form.

When $g=1$,
$$
F_1(x)=\Omega_1(y_0) -\frac12\log\Bigl(1- \Omega^{(2)}_0(y_0)\Bigr).
$$

When $g=2$,
\begin{eqnarray*}
F_2(x) &=& \Omega_2(y_0)-y_1^2y_0'-\frac{y_1^2}{2}+
      \frac{\Omega_1^{(2)}(y_0)}{2\bigl(1-\Omega_0^{(2)}(y_0)\bigr)}
      -\frac{1}{8\bigl(1-\Omega_0^{(2)}(y_0)\bigr)^{2}(1+y_0)^2} +
      \frac{1}{12(1+y_0)^2\bigl(1-\Omega_0^{(2)}(y_0)\bigr)} + \\ &+&
      \frac{1}{8(1+y_0)^2(1-\Omega_0^{(2)}(y_0)\bigr)^3}.
\end{eqnarray*}
}
\end{example}

\section{The Ordinary Euler Characteristic of
$\overline{\cal M}_g^n$}
\label{sec3}

In this section, we will express the ordinary Euler characteristic of
$\overline{\cal M}_g^n$ in terms of $\chi({\cal M}_g^n)$. This amounts
to the computation of $e(\Delta^o_G)$ for any stable graph $G$.  For
this purpose, we shall pursue previous work in \cite{Haz} and apply
some results of group cohomology theory.  In fact, $\Delta^o_G$ is a
rational $K(\Gamma(G),1)$; hence we have
$$
e(\Delta^o_G)=e(\Gamma(G)).
$$

Thus, the computation of $e(\Gamma(G))$ will follow from
a result in \cite{Bro}. Define a group $K$ to be
{\it geometrically WFL} if there is a contractible, finite,
dimensional, proper $K$-complex $Y$ such that there are only finitely
many cells of $Y$ under the action of $K$. Suppose, further, that $K$ 
has finitely many conjugacy classes of elements of finite order and 
for every element $\sigma$ in $K$ the centralizer $Z_K(\sigma)$ is 
geometrically WFL. Then the following holds.

\begin{theorem}(\cite{Bro})
\label{cobro}
For each $\sigma$ of finite order in $K$
$$
e(K)=\sum_{C_{\sigma}}\chi(Z_K(\sigma)),
$$
where the sum is over all conjugacy classes $C_{\sigma}$ of elements
of finite order in $K$, and $\chi(Z_K(\sigma))$ is the Euler characteristic of the group 
$Z_K(\sigma)$ in the sense of Wall (cf. \cite{B})
\end{theorem}

We shall apply Theorem \ref{cobro} to the group $\Gamma(G)$ for any stable graph $G$.
Let $Y_g^n$, $n \geq 1$, be the CW-complex introduced in
\cite{Har2}. $Y_g^n$ is a
contractible, finite dimensional complex such that the mapping class group
$\Gamma_g^n$ acts cellularly, with finite stabilizers and finitely many orbits.
For a graph $G$, consider the CW-complex given by the product
\begin{equation}
\label{complex}
Y(G):= \prod_{v \in G} Y_{g(v)}^{l(v)}.
\end{equation}

By the properties of $Y(G)$, the group $\Gamma(G)$ is geometrically WFL.
We shall prove that centralizers of elements of finite order in $\Gamma(G)$ are
geometrically WFL in Corollary \ref{mesac}. Since, as we shall see,
$\chi(Z_G(\sigma))$ can be computed in terms of the
characteristic of a group which is a finite extension of
products of mapping class groups, the ordinary Euler characteristic of the stratum
$\Delta^o_G$ is determined by $\chi({\cal M}_g^n)$.
Various algebraic manipulations will yield the final result.

\subsection{$e({\cal M}_g^{n+1})$}\label{sec4}

To exemplify the strategy above, we consider first the stable graph
$G$ with one vertex and $n+1$ legs. This will yield a
formula for the open locus of smooth pointed curves. In
this section we restrict our attention
to curves with at least two marked points. The remaining cases
are dealt with in \cite{Haz}.

Fix a genus $g$ topological oriented surface $S_{g,n+1}$ with $n+1$ marked 
points. Let $\sigma \in \Gamma_g^{n+1}$ be an element of finite order. As proved in 
\cite{Nie}, $\sigma$ may be represented by a periodic
homeomorphism $\mathfrak{}f$ of order $k$
which fix $p_i \in S_{g,n+1}$, $i=1, \ldots, n+1$.
Such a homeomorphism defines a branched covering
$$
\psi_{\mathfrak f}:S_{g,n+1} \rightarrow
H_{g,n+1}:=S_{g,n+1}/<{\mathfrak f}>,
$$
where $H_{g,n}$ has a natural structure of orbifold of
genus $h$. If
$p_1, \ldots, p_n, p_{n+1}, \ldots, p_{n+d+1}$ denote the ramification points, then
by Riemann-Hurwitz Formula we have:
\begin{equation}
\label{cond}
2g-1+n=k(2h-1+n+d) - \sum_i M_i,
\end{equation}
where $k \geq 1$, $h \geq 0$ and $1 \leq M_r < k, \, M_r|k$. The unramified covering
corresponding to $\psi_{\mathfrak f}$ is clearly
determined by a group homomorphism
\begin{equation}
\label{cond1}
\omega_{\sigma}: H_1(H_{g,n+1} - B) \rightarrow {\mathbb Z}/k{\mathbb Z},
\end{equation}
where $B$ is the branch locus of $\psi_{\mathfrak f}$. As a result, an
element of finite order in $\Gamma_g^{n+1}$ determines a homomorphism
$\omega_{\sigma}$ and integers $h,k,d, M_i$ satisfying
\eqref{cond}. On the other hand, it is easy to check that data
$\{h,k,d, M_i\}$ are sufficient to have an element of order $k$ in
$\Gamma_g^{n+1}$. Define, now, $\Gamma(H_{g,n+1})$ to be the group of
all isotopy classes of homeomorphisms of $H_{g,n+1}$, which fix the set
$\{p_1, \ldots, p_{n+1}\}$ and may permute $p_i$ and $p_j$ for $i,j
\geq n+2$ -- when they have the same monodromy. Set, further,
$$
\Gamma_g^{n+1}({\mathfrak f}):=\left\{\mathfrak{h}
\in \Gamma(H_{g,n+1}):
\omega_{\sigma} \circ \mathfrak{h} =
\omega_{\sigma}\right\}.
$$

Let $Z_{\sigma}$ and $N_{\sigma}$ be the centralizer and the normalizer of 
$\sigma$ in $\Gamma_g^{n+1}$, respectively. Analogously to Lemma 3 in \cite{Haz},
the following holds.
\begin{lemma}
\label{centro}
The groups $N_{\sigma}$ and $\Gamma_g^{n+1}({\mathfrak f})$ are related via the short 
exact sequence
$$
1\rightarrow {\mathbb Z}/k{\mathbb Z} \rightarrow N_{\sigma} \rightarrow \Gamma_g^{n+1}(f)\rightarrow 1.
$$
Moreover, the groups $N_{\sigma}$ and $Z_{\sigma}$ are geometrically WFL. 
\end{lemma}

By Lemma \ref{centro}, $\chi(Z_{\sigma})$ is well defined and can be
computed in terms of the Euler characteristic of $N_{\sigma}$. Similar
arguments to those in \cite{Haz} yield a closed formula for $e({\cal
M}_g^{n+1})$. Let us recall now some conventional notation. As
customary, we denote by $\phi$ and $\mu$ the Euler and the M\"{o}bius
arithmetic functions, respectively. Additionally, for any triple of
non-negative integers $k,l, \delta$ such that $l|k$ and $\delta|k$, we
set
\begin{equation}
\label{ecci}
c(k,l,\delta)=\frac{\phi(k/l)}{\phi\bigl(\delta/(\delta,l)\bigr)}\mu\bigl(\delta/(\delta,l)\bigr),
\end{equation}
where the number $(\delta,l)$ is the g.c.d. of $\delta$ and $l$. Then the following
holds.

\begin{theorem}
\label{statest}
For nonnegative integers $g,n$ such that $2g-1+n>0$, the ordinary
Euler characteristic of ${\cal M}_g^{n+1}$ is
\begin{equation}
e({\cal M}_g^{n+1})=\sum_{h,k,M_1, \ldots, M_d} \frac{\phi(k)}{k}
 \frac{\chi({\cal M}_h^{d+1+n})}{d!} k^{2h-1} 
 \sum_{\delta|k}\mu(\delta)\left(c(k,1,\delta)\right)^n \cdot
\prod_{r=1}^d c(k,M_r,\delta),
\label{formst}
\end{equation}
where $h,k,d, M_1, \ldots, M_d$ satisfy the following conditions:
$$
k \geq 1, \, \, h \geq 0, \, \, 1 \leq M_r < k, \, \, M_r|k,
$$
$$
2g-1+n=k(2h-1+n+d)-M_1-M_2- \ldots -M_d.
$$
\end{theorem}

\begin{remark}
Note that for $n=0$ Formula \eqref{formst} 
coincides with the one given in \cite{Haz}. 
\end{remark}

In Table 1 we give some values of $e({\cal M}_g^{n+1})$ for 
$3 \leq g \leq 10$ and $1 \leq n \leq 8$. In fact, all the values we get for 
$g=0,1,2$ coincide with the known ones. Finally, we further remark that 
Formula \eqref{formst} generates the same numbers as 
$\chi({\cal M}_g^{n+1})$ for $n \geq 2g+2$. This is consistent with 
the general fact that a smooth curve with at least $2g+3$ marked 
points is automorphism free; hence the two Euler characteristics coincide. 

\begin{center}
\begin{tabular}{r|r|r|r|r|r|r|r|r|}
$g \Big\backslash n$ & $1$ & $2$ & $3$ & $4$ & $5$ & $6$ & $7$ & $8$
 \\ \hline
 $3$ & $8$ & $6$ & $4$ & $-10$ & $30$ & $-660$ & $6540$ & $-79200$ \\ \hline
 $4$ & $-2$ & $-10$ & $-24$ & $-24$ &  $-360$ & $2352$ & $-37296$ & 
 $501984$  \\ \hline
 $5$ & $12$ & $26$ & $92$ & $182$ & $1674$ & $-16716$ 
 & $238980$ &  $-3961440$ 
 \\ \hline
 $6$ & $0$ & $-46$ & $-206$ & $188$ & $-7512$ & $124296$ 
 & $-2068392$ & $37108656$ 
 \\ \hline
 $7$ & $38$ & $120$ & $676$ & $-1862$  & $71866$ 
 & $-1058676$ & $21391644$ & $-422727360$ 
 \\ \hline
 $8$ & $-166$ & $-630$ & $-5362$ & $16108$ & $-680616$ 
 & $12234600$ & $-259464240$ & $5719946400$
 \\ \hline
 $9$ & $748$ & $2132$ & $29632$ & $-323546$ & $7462326$ & $-164522628$ & $3771668220$ & 
 $-90553767840$  \\ \hline
 $10$ & $-1994$ & $6078$ & $-213066$ & $4673496$ & $-106944744$ & $2559934440$ & $-64133209320$ &
  $1.6663E+12$ \\ \hline
\end{tabular}
\end{center}

\vskip 0.2cm
\centerline{\sc Table 1: Some values of $e({\cal M}_g^{n+1})$}

\subsection{The General Case}\label{sec5}

Analogously to Section \ref{sec4}, we shall give a formula for
$e(\overline{\cal M}_g^n)$. First of all, we arrange all
such
numbers in the generating function
\begin{equation}
\label{fform}
f(\lambda, y):= \sum_{\substack{g \geq 0 \\ n \geq 1 \\ 2g+n \geq 3}}
e(\overline{\cal M}_g^n)\lambda^{2g-2+n}
\frac{y^n}{n!}
\end{equation}
and express $f$ in terms of matrix integrals. To state the main formula,
we need some more notation. For any nonnegative integers $k$ and
$\delta|k$, define $V(k,\delta)$ to be the polynomial 
\begin{eqnarray}
\label{vform}
V(k,\delta)&=& c(k,1,\delta)\lambda^ky + \sum_{1 \leq m < k}
c(k,m,\delta) \lambda^{k-m}  + T(k,\delta)\lambda^k
\nonumber \\ &+& \sum_{\underset{r|k}{1 \leq r \leq k}}
c(k,r,\delta) x_r \lambda^k
\end{eqnarray}
in the variables $\lambda, y, x_1, x_2, \ldots, x_k$.  In
\eqref{vform}, note that
\begin{equation}
T(k,\delta)=
\left\{
\begin{array}{cl}
\frac{k}{2} & k\equiv 0 \, \, \hbox{mod}\, 2, \delta=1,2, \\
\\ & \\
0           & \hbox{otherwise},
\end{array}
\right.
\label{tform}
\end{equation}
and the other coefficients are defined in \eqref{ecci}.
Set, further,
\begin{equation}
\label{gform}
Q(\lambda,y,\underline{x})=\sum_{\delta|k}
\sum_{\underset{h,s}{2h+s \geq 3}}\phi(\delta)
\chi({\cal M}_h^s)(k\lambda^k)^{2h-2}\frac{V^s(k,\delta)}{s!}.
\end{equation}

Then the following holds.
\begin{theorem}
\label{matform}
\begin{equation}
\label{matform1}
f(\lambda,y)+ \sum_{g \geq 2} \left( \sum_{\substack{ h=1 \\
h|(g-1)}}^{g-1}e({\overline{\mathcal M}}_{h+1}) - e({\overline{\mathcal
M}}_g) \right) \lambda^{2g-2} = \log \left\{lim_{M \rightarrow \infty}
\frac{1}{(2\pi)^{M/2}}\int_{{\mathbb R}^M} \exp(Q)d\mu_M\right\},
\end{equation}
where
\begin{equation}
\label{misura}
d\mu_M= \exp\left( -\frac12 (x_1^2 + \hdots x^2_M)\right)
dx_1 \ldots dx_M.
\end{equation}
\end{theorem}

\begin{remark} Note that the results in Section \ref{sec4} can
be expressed in terms of generating functions, too. In fact,
set
$$
f_0(\lambda,y) = \sum_{\substack{ g,n \geq 0 \\ 2g+n \geq 3}}
\frac{e({\cal M}_g^n)}{n!}\lambda^{2g-2+n}y^{n}
$$
and
$$
v_0(k,\delta)= c(k,1,\delta)\lambda^ky +
\sum_{\underset{m|k} {1 \leq m < k}}c(k,m,\delta)\lambda^{k-m}.
$$

Then
$$ f_0(\lambda,y)- \sum_{g \geq 2} \left( \sum_{\substack{ h=1 \\
h|g-1}}^{g-1}e({\mathcal M}_{h+1}) - e({\mathcal
M}_g) \right) \lambda^{2g-2}= \sum_{\underset{h,s \geq 0}{2h+s \geq
3}} \sum_{1 \leq \delta|k}\phi(\delta) \chi({\cal
M}_h^s)(k\lambda^k)^{2h-2}\frac{V_0^s(k,\delta)}{s!}.
$$
\end{remark}

As sketched before, the proof of Theorem \ref{matform} is organized
as follows. To begin with, we will prove that centralizers
of elements of finite order in $\Gamma(G)$
are geometrically WFL - see Section \ref{sec6}. Henceforth,
we assume $G$ is a stable graph with at least two vertices.
Next, we will obtain a formula for
$e({\overline{\cal M}}_g^n)$ in terms of $\chi({\cal M}_h^r)$ for suitable values of $h$ and $r$ -
see Section \ref{sec6}. Finally, we will deduce Formula
\eqref{matform1} by standard techniques in matrix integral
theory.

\subsubsection{Centralizers of Elements of Finite Order in
$\Gamma(G)$}\label{sec00}

The elements of finite order in $\Gamma(G)$, and their centralizers,
are better understood if we describe the group $\Gamma(G)$
in an alternative way. For any stable graph $G$, let $S(G)$ be the surface
\begin{equation}
\label{essedi}
S(G)=\bigsqcup_{v \in V_{G}}S_{g(v), a(v)+b(v)},
\end{equation}
where  $S_{g(v), a(v)+b(v)}$ is defined in Section \ref{sec1}.
For nonnegative integers $\mathfrak{g},\mathfrak{b}$
denote by $V_{G}^{(\mathfrak{g},\mathfrak{b})}$ the
subset of vertices of $G$ such that $g(v)=\mathfrak{g}$,
$a(v)=0$ and $b(v)=\mathfrak{b}$. Any non-empty
$V_{G}^{(\mathfrak{g},\mathfrak{b})}$ is equipped with
a permutation action of the symmetric group
$\mathfrak{S}_{(\mathfrak{g},\mathfrak{b})}$. Define
$\Gamma(S(G))$ to be the semidirect product
$$ 
\widetilde{\Gamma} \rtimes
\prod_{V_{G}^{(\mathfrak{g},\mathfrak{b})} \neq \emptyset}
\mathfrak{S}_{(\mathfrak{g},\mathfrak{b})},
$$
where $\widetilde{\Gamma}$ consists of isotopy classes of
orientation-preserving diffeomorphisms of $S(G)$ which, for each
vertex $v$, fix the $a(v)$ points, but may permute the $b(v)$
points.

We now show that $\Gamma(G)$ is isomorphic to a subgroup $L$
of $\Gamma(S(G))$. Let $L$ be the set of elements of
$\Gamma(S(G))$
that are compatible with the
identifications of the marked points, which are
induced by $G$. In other words, if $q_i$ and $q_j$ are identified in
$S(G)$, then any $\mathfrak{h}$ in $L$ identifies
${\mathfrak h}(q_i)$ and ${\mathfrak h}(q_j)$. It is easy to check that
$\Gamma(G)$ is isomorphic to $L$. Clearly, any element
of $\Gamma(G)$ induces an element in $L$. Conversely, note that
${\mathfrak h}$ in $L$ satisfies
${\mathfrak h}(q_i)=q_{\theta(i)}$ and ${\mathfrak
h}(q_j)=q_{\theta(j)}$, where $\theta$ is a permutation of the set of
the $\sum_v b(v)$ marked points of $S(G)$. Accordingly,
${\mathfrak h}$ induces the element $(\widetilde{{\mathfrak h}}, \widetilde{\theta})$
in $\Gamma(G)$, where $\widetilde{\theta}$ is the automorphism of $G$
induced by $\theta$ and $\widetilde{{\mathfrak h}}=
\widetilde{\theta}\circ \mathfrak{h}$.

Since $\Gamma(G)$ can be viewed as a subgroup of $\Gamma(S(G))$,
an element $\sigma$ of finite order in $\Gamma(G)$ can
be realized (cf. \cite{Nie}) as a periodic diffeomorphism
of $S(G)$. The quotient of $S(G)$ by the group generated by $\sigma$
is a disconnected orbifold, $X$. Let us describe $X$ in detail.

$X$ has a finite number (say $p$) of connected components $X_i$ each
of which is an orbifold of genus $h_i, \, 1 \leq i \leq p$. Each $X_i$
has some natural marked points. The $n$ marked points of $S(G)$
corresponding to legs of $G$ are fixed by $\sigma$ and thus descend to
marked points of $X$.  Hence, each $X_i$ has $a_i \, \, (a_i \geq 0)$
points of this type, where $\sum a_i=n$. Another set of marked points
is given by the orbits of $\langle \sigma \rangle$ that contain points
of $S(G)$ corresponding to half-edges of $G$. Additionally, we mark
the points of $X_i$, which come from orbits of smaller length than
$|\langle \sigma \rangle|$. Each $X_i$ may have $d_i$ additional
marked points of this type.

This leads quite naturally to the construction of an {\it orbi-graph}
$H$ in the following way. The set of vertices $\{v_1, \ldots, v_p\}$
of $H$ has order $p$. Each $v_i$ corresponds to one of the
$X_i$'s. The edges of $H$ are recovered from the edges of $G$ if we
respect some compatibility conditions. More explicitly, let $\beta_i$
and $\beta_j$ be two half-edges of $H$. Suppose they correspond to
points $x_i$ and $x_j$ of $X$. If $y_i$ and $y_j$ are points of
$S(G)$ which correspond to $x_i$ and $x_j$, then
$\beta_i$ and $\beta_j$ are paired if and only if the
half-edges of $G$ corresponding to $y_i$ and $y_j$ are
paired too. Note that some edges of $G$ may
correspond to half-edges of
$H$. Thus, we denote by $b_i$ the number of marked
points in $X_i$ that correspond to the
half-edges issuing from $v_i$ and
giving rise to edges of $H$. On the other hand, we denote by $c_i$
the number of half-edges issuing from $v_i$ and different
from the $a_i+b_i+d_i$ half-edges listed so far. Finally, we denote by
$k_i$ the degree of the covering, say $\sigma_i$,
of $X_i$. An element of finite
order in $\Gamma(G)$ thus determines the following set of data:
\begin{equation}
\label{ec14}
p,\{k_1, \ldots, k_p\} ,\{h_1, \ldots, h_p\}, \{a_1, \ldots, a_p\},
\end{equation}
$$
\{b_1, \ldots, b_p\},\{c_1, \ldots, c_p\}, \{d_1, \ldots, d_p\}, \pi,
$$
where $\pi$ is the pairing among the $\sum_{i=1}^p b_i$
half-edges of $H$ induced by the compatibility conditions
described above.

Note that a stable graph $G$ and the elements of finite order
$\sigma_i$ determine group homomorphisms from $H_1(X_i^0)$ to the
cyclic group of order $k_i$, where $X_i^0$ is $X_i$ with all the
marked points removed. Moreover, these homomorphisms satisfy the
following conditions. In the sequel, when we refer to small loops
around a marked point $q$ of $X_i^0$ we mean an oriented loop (say
counterclockwise) around $q$, which is small enough so that it
encircles only the point $q$.
\begin{itemize}
\item [I)]\label{iuno} $(\omega_i(\alpha),k_i)=1$, for a small loop
$\alpha$ around each of the $a_i$ marked points;
\item [II)]\label{idue}
$\omega_i(\delta)\neq 0$, for a small loop around each of the $d_i$
marked points;
\item [III)]\label{itre}
$(\omega_i(\gamma),k_i) \equiv 0$ mod $2$, for a small loop $\gamma$
around each of the $c_i$ marked points.
\end{itemize}

Conversely, given a graph $H$ and data as in (\ref{ec14}), we would
like to recover a stable graph $G'$ of type $(g,n)$ for some $g$ and
$n$. Roughly speaking, we would like to reconstruct a topological
surface from $X$. For these purposes, we need some extra
information. First, we define group homomorphisms

\begin{equation}
\label{ec15}
\omega_i: H_1(X_i^0) \rightarrow \mathbb{Z}/{k_i\mathbb{Z}},
\end{equation} 
which satisfy conditions I), II) and III). A homomorphism
$\omega_i$ in (\ref{ec15}) determines a branched covering of $X_i$ of
degree $k_i$: the image of the loops around the marked points of $X$
determine the local monodromy of the covering. Let
\begin{equation}
\label{ec23}
N_j^i, \, \, 1 \leq j \leq b_i, \quad 
N_j^i|k_i, \, \, 1 \leq N_j^i \leq k_i,
\end{equation}
\begin{equation}
\label{ec25}
M_r^i, \, \, 1 \leq r \leq d_i, \quad 
M_r^i|k_i,\, \,  1 \leq M_r^i < k_i,
\end{equation}
\begin{equation}
\label{ec26}
R_s^i, \, \, 1 \leq s \leq c_i,
\end{equation}
be the number of points lying over each of the $b_i, d_i, c_i$ marked
points of $X_i$. Notice that each vertex of the graph $G'$ has
\begin{equation}
\label{he}
\eta_i = a_i + \sum_{j=1}^{b_i} N_j^i + \sum_{s=1}^{c_i} R_s^i
\end{equation}
half-edges. Note that the sums in (\ref{he}) should be considered
empty when $b_i \, (c_i \, \hbox{or} \, d_i)$ are zero. Next, we
assign a pairing $\widetilde{\pi}$ among the points in the orbit of
each of the $b_i$ points.  Clearly, this pairing has to be compatible
with the pairing $\pi$ of the graph $H$.  Indeed, suppose that
$\beta_j^i$ and $\beta_{j'}^{i'}$ are paired via $\pi$. By abuse of
notation, we still denote the marked points corresponding to them by
$\beta_j^i$ and $\beta_{j'}^{i'}$.  If there are $N_j^i$ (resp.
$N_{j'}^{i'}$) points lying above $\beta_j^i$
(resp. $\beta_{j'}^{i'}$), then $N_j^i=N_{j'}^{i'}$.

We can now associate a graph $G'$ to the covering of $X$ determined by
the homomorphisms $\omega_i$. The number, $v(G')$, of vertices of $G'$
is $\sum_{i=1}^pP_i$, where $P_i$ is the number of components lying
over $X_i$. The total number, $e(G')$, of edges is determined as
follows. There is only one way to lift the $c_i$ half-edges of $H$, so
$G'$ has $\sum_{i,s}R_s^i/2$ edges of this type. Let $e$ be an edge
obtained by pairing two half-edges $e_{-}$ and $e_{+}$. Denote by
$N(e_{-})$ and $N(e_+)$ the number of points over (the marked points
corresponding to) $e_{-}$ and $e_{+}$. By the discussion above, there
exist $N(e_{-})=N(e_{+})$ possible pairings $\widetilde{\pi}$ for
$G'$. As a consequence, there are $1/2\sum_{i,j}N^i_j$ edges of this
type. We observe that each vertex of $G'$ corresponds to a genus $g_i$
covering of $X_i$, where $g_i$ is determined by the Riemann-Hurwitz
formula, namely:

\begin{equation}
\label{eqt1}
P_i(2-2g_i)=k_i(2-2h_i -a_i -b_i -c_i -d_i)+a_i +  \sum_{j=1}^{b_i} N_j^i +
\sum_{r=1}^{d_i} M_r^i + \sum_{s=1}^{c_i} R_s^i.
\end{equation}

If we set

\begin{equation}
\label{ec20}
n :=\sum_{i=1}^p a_i, \quad g:=\sum_{i=1}^p P_ig_i +e(G') -v(G')+1,
\end{equation}
we get

\begin{eqnarray}
\label{ec16}
2g-2+n &=&\left(\sum_{i=1}^p k_i\Bigl(2h_i -2 +a_i +b_i +c_i
       +d_i\Bigr)\right)-\sum_{i,r}M^i_r.
\end{eqnarray}
 
The graph $G'$ does not have a total ordering of all
the $a_1+ \ldots +a_p$ points. Thus, we give a total ordering
${\bf{O}}$ to all these points and there are clearly $n!$ of these.
We call the result $G$, and the following now holds.

\begin{proposition}
\label{ec17}
The graph $G$ is stable if and only if
\begin{equation}
\label{newstab}
k_i(2h_i -2 + a_i +b_i + c_i + d_i) - \sum M^i_r > 0
\end{equation}
holds for all $1 \leq i \leq p$.  Furthermore, this condition is equivalent to the stability of $H$ except for the case where $h_i=a_i=c_i=0$, $b_i=1$, $d_i=2$, $k_i$ is even and $M^1_i=M_i^2=k_i/2$.  If these values hold for some $i$, then $G$ is not stable, even though $H$ is.

\end{proposition}

\proof Let $a_i+b_i+c_i+d_i$ be the number of half-edges at the $i$-th
vertex of $H$ and let $\eta_i$ be the number of half-edges of $G\prime$ as in
(\ref{he}). Suppose first that \eqref{newstab} holds.  This means that $2h_i-2
+a_i+b_i+c_i+d_i >0$ for $1 \leq i \leq p$. By (\ref{eqt1}) we have

$$ P_i(2g_i-2+\eta_i) \geq P_i(2g_i -2)+\eta_i =k_i(2h_i-2
+a_i+b_i+c_i+d_i) - \sum M^i_r >0,
$$
so $G$ is stable. 

Conversely, suppose that $H$ is not stable at a vertex $v_i, 1 \leq i \leq
p$. Then either
$$
h_i=1, \, \, a_i+b_i+c_i+d_i=0
$$ 
or
$$
h_i=0, \, \, a_i+b_i+c_i+d_i \leq 2.
$$ 

Let $\widetilde{X}_i$ be one of the $P_i$ components lying above
$X_i$. In the former case, there are no branch points, so
$\widetilde{X}_i$ must be a torus with no marked points implying that
$G$ is also not stable. In the latter case, $\widetilde{X}_i$ must be a
cover of $X_i$ branched at most over two points and again $G$ cannot
be stable.

\qed

 In other words, a pair $(G,\sigma)$ of a stable graph and an element
of finite order in $\Gamma(G)$ is determined by i) a stable graph $H$
that satisfies \eqref{newstab}, ii) a collection of data as in
(\ref{ec14}), iii) group homomorphisms $\omega_i$ satisfying
conditions I), II), III), and finally iv) a total ordering $\bf{O}$
and a pairing $\pi$.

This said, it is now easier to study properties of the
centralizer $Z_G(\sigma)$ of an element $\sigma$
of finite order in $\Gamma(G)$.

Let $N_G(\sigma)$ be the normalizer  of $\sigma$ in
$\Gamma(G)$. Denote by $H$ the graph associated with $(G, \sigma)$.
Since $H$ is a stable, we define the group $\Gamma(H)$ as in (\ref{arieccote}), i.e.
\begin{equation}
\label{ec18}
\prod_{i=1}^p\Gamma_{h_i, a_i+b_i+c_i+d_i} \rtimes Aut(H), 
\end{equation}
where $Aut(H)$ is the group of automorphisms of $H$ generated by
\begin{itemize}
\item automorphisms of $H$ which may permute two vertices $v_i$ and
$v_j$ only if $h_i=h_j, a_i=a_j=0, b_i=b_j, c_i=c_j$ and $d_i=d_j$,
\item permutations of the $d_i$ marked points for each $i, 1 \leq i \leq p$.
\end{itemize}

Any element $\mathfrak{h}$ in $\Gamma(H)$ induces an automorphism of
the fundamental group of $X$ which we denote by
$\mathfrak{h}_{*}$. Thus, we set
\begin{equation}
\label{stabilizza}
\Gamma(H,\omega_{\sigma}): =\{\mathfrak{h} \in \Gamma(H): \omega_{\sigma} \circ \mathfrak{h}_*=\omega_{\sigma}\},
\end{equation}
where $\omega_{\sigma}$ is any of the group homomorphisms determined by
$\sigma$.   

\begin{proposition}
The groups $N_G(\sigma)$ and $\Gamma(H, \omega)$ are
related via the short exact sequence of groups
\begin{equation}
\label{exact}
1\rightarrow C_k \overset{t_1}{\rightarrow} N_G({\sigma}) \overset{t_2}{\rightarrow}\Gamma(H,\omega_{\sigma}) \rightarrow 1.
\end{equation}
\label{mesa}
\end{proposition}

\proof The group $C_k$ is the cyclic group of order $k$ generated by $\sigma$
and $t_1$ is the inclusion homomorphism. The map $t_2$ is defined as
follows. Pick an element ${\mathfrak h}$ in $N_{G}(\sigma)$. By
the same arguments adopted in \cite{Haz} (Lemma 3), there exists a
preserving-orientation diffeomorphism
${\mathfrak f}_{{\mathfrak h}}$ of $S(G)$
such that

\begin{equation}
\label{coniugio}
{\mathfrak f}_{{\mathfrak h}}\sigma {\mathfrak f}_
{{\mathfrak h}}^{-1} =\sigma^j \, \,
\hbox{for some $j$}.
\end{equation}
  
Then $t_2({\mathfrak h})$ is defined as the isotopy
class of ${\mathfrak f}_{{\mathfrak
h}}$. By condition (\ref{coniugio}) and the
definition of $\Gamma(H,\omega_{\sigma})$,
${\mathfrak f}_{{\mathfrak h}}$ yields an
element in $\Gamma(H,\omega_{\sigma})$. 

\qed

\begin{corollary}
\label{mesac}
Centralizers of elements of finite order in $\Gamma(G)$ are geometrically WFL.
\end{corollary}

\proof The group $\Gamma(H,\omega_{\sigma})$ is a finite index
subgroup of $\Gamma(H)$ so it acts on the CW-complex $Y(G)$ introduced
in (\ref{complex}). The exact sequence in (\ref{exact})
gives an action of $N_G(\sigma)$ and $Z_G(\sigma)$ on
$Y(G)$; so they are both WFL.

\qed

\subsubsection{A Formula For
$e({\overline{\cal M}}_g^n)$} \label{sec6}

In this section we give a formula for
$e({\overline{\cal M}}_g^n)$ when $n \geq 1$. This will be used
in the next section for the generating function $f$.
In what follows, we adopt the same notation as that in
Section \ref{sec00}.

By Theorem \ref{cobro}, we have
\begin{equation}
\label{eform}
e(\mgn)=\sum_G e(\Delta^o_G)=\sum_G \sum_{C_\sigma}
\chi(Z_{G}(\sigma)).
\end{equation}

Since each $(G, \sigma)$ determines an orbi-graph $H$, we rewrite
$e(\mgn)$ as a sum over pairs $(H,\omega_i)$, where $H$ is an
orbi-graph and $\omega_i$ are the group homomorphisms introduced in
(\ref{ec15}). For each $k_i$ we consider the action of $\Gamma(H)$ on
the set $A_{H}(k_i)$ of all homomorphisms $\omega_i$ in (\ref{ec15}),
which satisfy conditions I), II) and III). In particular, if $\omega_i
\in A_{H}(k_i)$, we denote by $\Gamma(H,\omega_i)$ the stabilizer of
$\omega_i$ under the action of $\Gamma(H)$. Notice that when
$\omega_i$ is induced by $\sigma \in \Gamma(G)$, $\Gamma(H,\omega_i)$
is the group in (\ref{stabilizza}).

Fix an orbi-graph $H$ and data $\bf{O}$ and $\pi$ described in
Section \ref{sec5}. Let $\Lambda_i$ be a set of representatives of the
conjugacy classes $C_{\sigma_i}$, and suppose that the quotient of
$S(G)$ by the $\sigma_i$'s is isomorphic to $X$. Thus
$$
e(\mgn)=\sum_p \sum_{k_1, \ldots, k_p} \, \sum_{C_{\sigma_1}, \ldots, C_{\sigma_p}} \, \sum_{C_{\rho_1}, \ldots C_{\rho_p}} \prod_{i=1}^p \chi(Z_G(\sigma_i)),
$$
where $S_{\sigma_i}$ is a set of representatives of the orbits
of
$$
G_{\sigma_i}=\{\sigma_i^n: (n,k)=1\}
$$
under the action of $N_G(\sigma_i)$ by conjugation. If
${\mathcal O}_j$ is such an orbit, then we get
\begin{equation}
\label{hofattoluna}
\sum_{\rho \in S_{\sigma_i}}\chi\bigl(Z_{G}(\rho)\bigr)=
\sum_{\rho \in S_{\sigma}}|\mathcal{O}_j|\chi\bigl( N_{G}
(\rho)\bigr)=\phi(k)\chi\bigl( N_{G}(\sigma)\bigr).
\end{equation}

Hence we have:
$$              
\sum_p \sum_{k_1, \ldots, k_p} \sum_{C_{\sigma_1}, \ldots,
C_{\sigma_p}}\prod_{i=1}^p\frac{\phi(k_i)}{k_i}\chi(\Gamma(H,
\omega_{\sigma_i})).
$$
By the short exact sequence 
$$
1 \rightarrow \prod_{i=1}^p\Gamma_{h_i, a_i+b_i+c_i+d_i}  \rightarrow \Gamma(H) \rightarrow
Aut(H) \rightarrow 1,
$$
we get
\begin{equation}
\label{bohhhh}
 \chi\Bigl(\Gamma(H,\omega_{\sigma_i})\Bigr)=
\bigl[ \Gamma(H): \Gamma(H,\omega_{\sigma_i})\bigr]
\dfrac{\prod_{i=1}^p\Gamma_{h_i, a_i+b_i+c_i+d_i}} {|Aut(H)|}=
|{\mathcal{O}}(\omega_{\sigma_i})|
\dfrac{\prod_{i=1}^p\Gamma_{h_i,
a_i+b_i+c_i+d_i}} {|Aut(H)|},
\end{equation}
where ${\mathcal{O}}(\omega_{\sigma_i})$ is the orbit of
$\omega_{\sigma_i}$ under the action of $\Gamma(H)$ on
$A_{H}(k_i)$.

The expression in (\ref{bohhhh}) can be further
simplified. We say that
$H$ is ordered if it has an ordering on its vertices and an ordering
on each collection of the half-edges going out of a vertex $v$.  Let
${\bf H'}$ be the set of all ordered (disconnected)
stable orbi-graphs. The group
$$
\mathfrak{S}_p \ltimes \prod_{i=1}^p\bigl(\mathfrak{S}_{a_i} \times
\mathfrak{S}_{b_i} \times \mathfrak{S}_{c_i} \times \mathfrak{S}_{d_i}
\bigr)
$$ acts on ${\bf H}'$ with orbits equal to orbi-graphs. To determine
the ordered orbi-graph $H$, we only need to enumerate
$$
\begin{array}{ll}
(h_1, \ldots, h_p), \, \, h_i \geq 0, & (k_1, \ldots, k_p), \, \, k_i \geq 1, \\
(a_1, \ldots, a_p), \, a_i \geq 0, &
(b_1, \ldots, b_p), \, \, b_i \geq 0, \,\, \sum_{i=1}^pb_i \equiv 0 \, \,\hbox{mod $2$}, \\
(c_1, \ldots, c_p), \, \,c_i \geq 0, \, \, c_i \geq 0, &
(d_1, \ldots, d_p), \, \, d_i \geq 0,
\end{array}
$$
together with a pairing $\pi$ of the numbers $1, \ldots, b_1 + \ldots
+b_p$. Denote by $\underline{N}$ the vector
$$
(N^1_1, \ldots, N^1_{b_1}, \ldots, N^p_1, \ldots, N^1_{b_p}),
$$
where $N^i_j$ is defined in (\ref{ec23}). By (\ref{ec16}),
(\ref{ec20}), (\ref{ec25}), and (\ref{ec26}), we get
\begin{equation}
e(\mgn)= \sum_p \sum_{k_1, \ldots, k_p} \frac{1}{p!}\sum_{\omega_i \in A_H(k_i)}\frac{\phi(k_i)}{k_i} \sum_{\substack{h_i,a_i,b_i,c_i,d_i \geq 0 \\
              b_1 +\ldots + b_p \, even \\ 
              2h_i + a_i + b_i +c_i +d_i \geq 3}}
\prod_{i=1}^p\frac{\chi({\cal M}_{h_i}^{a_i+b_i+c_i+d_i})}{a_i!b_i!c_i!d_i!}.
\label{ec27}
\end{equation}

We can simplify the sum in (\ref{ec27}) if we
enumerate all the elements in $A_H(k_i)$. A homomorphism $\omega_i$ from
$H_1(X_i^0) \rightarrow
\mathbb{Z}/k_i\mathbb{Z}$ is determined by assigning its values on a
basis of $H_1(X_i)$ and on each small oriented loop around any of the
marked points in $X_i$. Therefore, for any $k_i$, the number of such
homomorphisms can be computed as follows. The number of values that
$\omega_i$ can assume on a basis of $H_1(X_i)$ is $\prod_i^p
k_i^{2h_i}$. For the values on the loops around the marked points we
introduce the following notation. Denote by $<\alpha_l^i>,
<\beta_j^i> <\gamma_s^i>, <\delta_r^i>$ the loops around the $a_i,
b_i, c_i, d_i$ marked points. Then $\omega$ is
determined by assigning the following elements of
$\mathbb{Z}/k_i\mathbb{Z}$:

\begin{enumerate}
\item $\omega(<\alpha_l^i>)=A_l^i$, with $(A_l^i,k_i)=1, \
,\, \, 1 \leq l \leq a_i$,
\item $\omega(<\beta_j^i>)=B^i_j$, with $(B_j^i, k_i)=N^i_j, \, \,
N_j^i|k_i, \, \,1 \leq N_j^i \leq k, \, \, 1 \leq j \leq b_i$,
\item $\omega(<\gamma_s^i>)=C_s^i$, with $(C_s^i,k_i) \equiv 0$ mod
$2$, $\, \, 1 \leq s \leq c_i$,
\item $\omega(<\delta_r^i>)=D_r^j$, with $(D_r^i,k)=M^i_r, \, \,
M^i_r|k, \, \, 1 \leq M^i_r < k, \, \, 1 \leq r
\leq d_i$.
\end{enumerate}

These requirements depend on the conditions I), II) and  III) satisfied by
$\omega_i$. Moreover, the image under $\omega_i$ of the relation among
cycles in $H_1(X_i)$ yields the additional constraint
$$
\sum_{i,l}A^i_l + \sum_{j,i}B^i_j + \sum_{s,i}C^i_s + \sum_{r,i}D^i_r
\equiv 0 \, \hbox{mod $k_i$}.
$$

Define
${\bf T}\bigl(\{a_i\}, \{b_i\}, \{c_i\}, \{d_i\}, \{N_j^i\},
\{M_r^i\}\bigr)$ to be the order of the set
$$
\Bigl\{(A_1^1,A_2^1, \ldots, D_{d_p}^p) : A_1^1=1,\, \,  (A_l^i,k_i)=1, \, \, (B^i_j,k)=N_j^i,
$$
$$
(C^i_s,k_i) \equiv 0 \, \, \hbox{mod $2$}, \, \, (D^i_r,k_i)=M^i_r,
$$
$$
\sum_{i,l}A^i_l + \sum_{j,i}B^i_j + \sum_{s,i}C^i_s + \sum_{r,i}D^i_r
\equiv 0 \, \, \hbox{mod $k_i$}\Bigr\}.
$$

Accordingly, if $a_i \geq 1$, the number of homomorphisms $\omega_i$ is
\begin{equation}
\label{ec29}
\prod_{i=1}^p(k_i^{2h_i})\phi(k_i){\bf T}
\bigl(\{a_i\}, \{b_i\}, \{c_i\}, \{d_i\}, \{N_j^i\}, \{M_r^i\}\bigr).
\end{equation}

By standard facts in elementary number theory (see \cite{Haz}) we have
$$
\phi(k_i){\bf T}\bigl(\{a_i\}, \{b_i\}, \{c_i\}, \{d_i\}, \{N_j^i\},
\{M_r^i\}\bigr)=
$$

\begin{equation*}
\prod_{i=1}^p\frac{1}{k_i}\sum_{\substack{\zeta \\ \zeta^k=1}}
\sum_{\substack{0 \leq A_l^i <k_i \\ (A_l^i,k_i)=1}} \, 
\sum_{\substack{0 \leq B_j^i <k_i \\ (B_j^i,k_i)=N_j^i}}
\sum_{\substack{0 \leq C_s^i <k_i \\ (C_s^i,k) \equiv 0 \, \, \hbox{mod $2$}}} \,
\sum_{\substack{0 \leq D_r^i <k_i \\ (D_r^i,k)=M_r^i}}
\zeta^{\sum_l A_l^i+ \sum_j B_j^i+ \sum_s C_s^i +\sum_r D_r^i}=
\end{equation*}
\begin{equation*}
\prod_{i=1}^p \frac{1}{k_i}\sum_{\substack{\zeta \\ \zeta^k_i=1}} \prod_{l=1}^{a_i}\left(
\sum_{\substack{0 \leq s < k_i \\ (s,k_i)=1}}\zeta^s\right) \prod_{j=1}^{b_i} 
\left(\sum_{\substack{0 \leq s < k_i \\ (s,k_i)=N_j^i}} \zeta^s\right)
\prod_{l=1}^{c_i}\left(
\sum_{\substack{0 \leq s < k_i \\ (s,k_i) \equiv 0 \, \, \hbox{mod $2$}}} \zeta^s
\right) 
\prod_{j=1}^{d_i} 
\left(\sum_{\substack{0 \leq s < k_i \\ (s,k_i)=M_r^i}}\zeta^s \right).
\end{equation*}

\bigskip

\begin{lemma}
\label{mipiaceproprio}
\begin{itemize}
\item [i)] If $k$ is even,
$$
\sum_{\substack{(r,k) \equiv 0 \, \, \hbox{mod $2$} \\ 0
\leq r <k}}\zeta^r = \left\{
\begin{array}{ll}
k/2 & \hbox{for} \, \, \zeta=1,-1, \\
0 & \hbox{otherwise}.
\end{array}
\right.
$$

\item [ii)] For any pair $l,\delta$ of divisors of $k$
and $\zeta$ a primitive $\delta$-th root of unity, we have
$$
\sum_{\substack{(r,k)=l \\ 0 \leq r <k}}\zeta^r=c(k,l,\delta),
$$
where $c(k,l,\delta)$ is defined in (\ref{ecci}).
\end{itemize}
\end{lemma}

\proof i) Since $k$ is even, we have
$$
\sum_{\substack{(r,k) \equiv 0 \, \, \hbox{mod $2$} \\ 0 \leq r <k}}\zeta^r = 
1+\zeta^2 +\zeta^4 + \ldots + \zeta^{k-2},
$$
which is zero unless $\zeta=1$ or $-1$ in which case it equals $k/2$. 

ii) It follows from the definition of the Mobi\"{u}s function.

\qed

By Proposition (\ref{mipiaceproprio}) i) and ii), we have:

\begin{equation}
\phi(k_i){\bf T}\bigl(\{a_i\}, \{b_i\}, \{c_i\}, \{d_i\}, \{N_j^i\},
\{M_r^i\}\bigr)= \prod_{i=1}^p \frac{1}{k_i}
\sum_{\substack{ \delta |k_i \\ 1 \leq \delta < k_i}}
c(k_i,1,\delta)^{a_i}\cdot  
\end{equation}
\begin{equation}
\cdot \prod_{j=1}^{b_i} c(k_i,N_j^i,\delta)\gamma(k_i,
\delta,c_i)\prod_{r=1}^{d_i}c(k_i,M_r^i,\delta),
\label{ec30}
\end{equation}
where 
\begin{equation}
\label{malena}
\gamma(k, \delta, c)= \left\{
\begin{array}{lcl}
\phi(\delta) & & c=0, \\
& & \\
0 & & k \equiv 1 \, \hbox{mod} \, 2, \, c>0, \\
& & \\
0 & & k \equiv 0 \, \hbox{mod} \, 2, \, c>0, \, \delta >2, \\
& & \\
(k/2)^c & &  k \equiv 0 \, \hbox{mod} \, 2, \, c>0, \delta=1,2.
\end{array}
\right.
\end{equation}

As a result, the following holds.
\begin{theorem}
\label{esoeso}
For any integers $g \geq 0$ and $n \geq 1$ such that $n > 2-2g$,
the ordinary Euler characteristic of $\mgn$ is given by
$$
e(\overline{\cal M}_g^n)= n! \sum_{p=1}^{2g-2+n} \frac{1}{p!}
\sum_{\underset{b_1 + \ldots +b_p \, even}
{b_1, \ldots, b_p}} \sum_{\underline{N}}
\sum_{\pi} \prod_{e} N(e)\prod_i
\sum_{\underset{2h_i+a_i+b_i+c_i+d_i \geq 3}
{h_i, a_i, d_i}} \sum_{k_i, M_i^r}  
\frac{\chi({\cal M}_{h_i}^{a_i+b_i+d_i+c_i})}
{d_i! a_i! b_i!c_i!} \cdot
$$
$$\cdot k_i^{2h_i-2} \sum_{\delta_i|k_i}
\Bigl(c(k_i,1,\delta_i)^{a_i}\Bigr)\prod_{j=1}^{b_i}
c(k_i,N_i^j,\delta_i) \gamma(k_i,\delta_i,c_i)
\prod_{r=1}^{d_i}c(k_i,M_i^r,\delta_i)\Bigr), 
$$
where
$$
2g-2+n=\sum_{i=1}^p\left(k_i(2h_i-2+a_i+b_i+c_i+d_i)\right)-\sum_{i=1}^p(M_i^1+ \ldots + M_i^{d_i});
$$
$$                                 
M_i^r|k_i,  \, \, M_i^r<k_i, \, \, k_i\leq 1;
$$
$$                                  
0 \leq h_i \leq g, \, \, \sum_{i=1}^p a_i=n;
$$
$$                                   
1 \leq N_i^j \leq k_i, \, \, N_i^j|k_i, \, \,
N_i^j=N^{j'}_{i'};
$$
$$
\pi \, \hbox{a connected pairing of the numbers} \,\,\, 1, \ldots,
b_1 +\hdots b_p, \quad  b_i \geq 1.
$$
\end{theorem}                             

\subsubsection{The Proof of Theorem \ref{matform}}
\label{sec7}
This section is devoted to proving Theorem \ref{matform}. Basically,
we shall apply Wick's Lemma to deduce \eqref{matform1}. 

By \eqref{ec16} and Theorem \ref{esoeso}, the generating series
$$
\sum_{\substack{ g \geq 1, n \geq 0 \\ n \geq 2g+2}}e(\overline{\cal
M}_g^n)\lambda^{2g-2}\frac{y^n}{n!}
$$
is equal to
\begin{equation}
\label{fine}
\sum_{p \geq 1} \frac{1}{p!}\sum_{k_1, \ldots, k_p}
\sum_{\substack{b_1, \ldots, b_p \\ b_1 + \hdots b_p \, \,
\hbox{even}}} \sum_{\underline{N}, \pi} \prod_{e} N(e)\cdot
\end{equation}
\begin{equation}
\cdot \prod_i \sum_{\underset{2h_i+a_i+b_i+c_i+d_i \geq 3}
{h_i, a_i, d_i}} \sum_{k_i, M_i^r} \frac{\chi({\cal M}_{h_i,
a_i+b_i+d_i+c_i})} {d_i! a_i! b_i!c_i!}(k_i\lambda^{k_i})^{2h_i-2}
\sum_{\delta_i|k_i}
\Bigl(\Bigl(c(k_i,1,\delta_i)\lambda^{k_i}y\Bigr)^{a_i}
\Bigr) \cdot
\end{equation}
\begin{equation}
\cdot \prod_{j=1}^{b_i}\Bigl(c(k_i,N_i^j,\delta_i)
\lambda^{k_i}\Bigr)
\Bigl(\gamma(k_i,\delta_i,c_i)\lambda^{k_i}\Bigr)
\prod_{r=1}^{d_i}\Bigl(c(k_i,M_i^r,\delta_i)
\lambda^{k_i - M^i_r}\Bigr)\Bigr).
\end{equation}

Note that
$$
\gamma(k_i,\delta_i,c_i)\lambda^{k_i}= \phi(\delta_i)
\Bigl(T(k_i, \delta_i)\lambda^{k_i}\Bigr)^{c_i},
$$ where $T(k_i, \delta_i)$ is defined in \eqref{tform}. Moreover, all
the indices $b_i$ in \eqref{fine} are positive integers. Looking at
the expansion above, we define the generating series
\begin{equation}
\widehat{f}(\lambda,y)= \sum_{p \geq 1} \frac{1}{p!}\sum_{k_1, \ldots,
k_p} \sum_{\substack{b_1, \ldots, b_p \geq 0\\ b_1 + \hdots b_p \, \,
\hbox{even}}} \sum_{\underline{N}, \pi} \prod_{e} N(e)\cdot
\label{fine1}
\end{equation}
\begin{equation*}
\cdot \prod_i \sum_{\underset{2h_i+a_i+b_i+c_i+d_i \geq 3}
{h_i, a_i, d_i}} \sum_{k_i, M_i^r} \frac{\chi({\cal M}_{h_i,
a_i+b_i+d_i+c_i})} {d_i! a_i! b_i!c_i!}(k_i\lambda^{k_i})^{2h_i-2}
\sum_{\delta_i|k_i}
\Bigl(\Bigl(c(k_i,1,\delta_i)\lambda^{k_i}y\Bigr)^{a_i}
\Bigr) \cdot
\end{equation*}
\begin{equation*}
\cdot \prod_{j=1}^{b_i}\Bigl(c(k_i,N_i^j,\delta_i)
\lambda^{k_i}\Bigr)
\Bigl(\gamma(k_i,\delta_i,c_i)\lambda^{k_i}\Bigr)
\prod_{r=1}^{d_i}\Bigl(c(k_i,M_i^r,\delta_i)
\lambda^{k_i - M^i_r}\Bigr)\Bigr).
\end{equation*}

Clearly, we have 
$$ \widehat{f}(\lambda,y)= \sum_{g \geq 2} u_g \lambda^{2g-2}+
\sum_{\substack{ g \geq 1, n \geq 0 \\ n \geq 2g+2}}e(\overline{\cal
M}_g^n)\lambda^{2g-2}\frac{y^n}{n!}
$$
whence
$$ f(\lambda,y) = \widehat{f}(\lambda,y)- \sum_{g \geq 2} \left( u_g -
e({\overline{\mathcal M}}_g) \right) \lambda^{2g-2}.
$$

Theorem \eqref{matform} will be completely proved if we show that i)
$\widehat{f}(\lambda,y)$ is equal to the right hand side of
\eqref{matform1} and ii) that the following holds:
\begin{equation}
\label{fine3}
u_g = \sum_{\substack{h=1 \\ h|(g-1)}}^{g-1} e({\overline{\mathcal M}}_{h+1}).
\end{equation}

To prove i) we argue as follows. Set
$$
\widehat{Q}(k,b,N_1, \ldots, N_b)= \sum_{\delta | k}\phi(\delta)
\prod_{j=1}^b \sqrt{N_j}c(k,N_j,\delta)\lambda^k \left(
\sum_{\underset{h,s}{2h+s+b \geq 3}} \chi({\cal M}_h^{b+s})
(k\lambda^k)^{2h-2} \cdot \right.
$$
$$
\cdot \left. \Bigl( c(k,1,\delta)\lambda^ky +T(k,\delta)
\lambda^k +
\sum_{\underset {m|k}{1 \leq m < k}}
c(k,m,\delta)\lambda^{k-m}\Bigr)^s \right).
$$

Then we have
\begin{equation}
\label{fine0}
\widehat{f}(\lambda, y)= \sum_{p \geq 1}
\frac{1}{p!}\sum_{k_1, \ldots, k_p}
\sum_{\substack{b_1, \ldots, b_p \geq 0 \\
b_1 + \hdots b_p \, \, \hbox{even}}} \, \, 
\sum_{\underline{N}, \pi}\, \, 
\prod_{i=1}^p \widehat{Q}(k_i, b_i, N^1_i, \ldots, N^p_i). 
\end{equation}

If we now expand $\exp(\widehat{f})$, we get
\begin{equation}
\label{fineform}
\sum_{q \geq 0} \frac{1}{q!}
\sum_{k_1, \ldots, k_q \geq 0}
\sum_{\substack{N^j_i|k_i \\ b_1, \ldots, b_q \geq 0}}
\sum_{ \pi} 
\prod_{i=1}^q \widehat{Q}(k_i, b_i, N^1_i, \ldots,
N^{b_i}_i).,
\end{equation}
where $\pi$ is a not necessarily connected pairing.
For any positive integer $r$ denote by $N_r$ the number of
those $N^j_i$'s that equal $r$. Set also
$$
{\mathcal N}_r= \left\{
\begin{array}{cl}
(N_r-1)!! & N_r \, \, \, \hbox{even} \\
0 & \hbox{else}.
\end{array}
\right.
$$

Formula \eqref{fineform} can be written as
$$
\sum_{q \geq 0} \frac{1}{q!}
\sum_{k_1, \ldots, k_q \geq 0}
\sum_{\substack{N^j_i|k_i \\ b_1, \ldots, b_q \geq 0}}
\, \prod_{r \geq 1}{\mathcal N}_r 
\, \prod_{i=1}^q \widehat{Q}(k_i, b_i, N^1_i,
\ldots, N^{b_i}_i).
$$

For any integer $k, b, N_1, \ldots, N_b, s$
define $R(\lambda, y, x_{N_1}, \ldots, x_{N_b})$
to be the polynomial
$$
R(\lambda,y, x_{N_1}, \ldots, x_{N_b})=
 \sum_{\delta | k}\phi(\delta)
\prod_{j=1}^b \sqrt{N_j}c(k,N_j,\delta)x_{N_j}
\lambda^k \left(
\sum_{\underset{h,s}{2h+s+b \geq 3}} \chi({\cal M}_h^{b+s})
(k\lambda^k)^{2h-2} \cdot \right.
$$
$$
\cdot \left. \Bigl( c(k,1,\delta)\lambda^ky +T(k,\delta)
\lambda^k +\sum_{\underset {m|k}{1 \leq m < k}}
c(k,m,\delta)\lambda^{k-m}\Bigr)^s \right).
$$

By Wick's Lemma, we  get
\begin{equation}
\exp(\widehat{f})=
lim_{M \rightarrow \infty}
\frac{1}{\sqrt{(2\pi)^M}} \int_{{\mathbb R}^M}
\exp \Bigl(
\,\sum_{\underset{k \geq 1}{b \geq 0}}
\,\sum_{N_i|k} 
R(\lambda,y, x_{N_1}, \ldots, x_{N_b})\Bigr)\,d\mu_M,
\end{equation}
where $d\mu_M$ is the gaussian measure as
in \eqref{misura}. It is an easy exercise
(left to the reader) to check that
$$
\sum_{\underset{k \geq 1}{b \geq 0}}
\sum_{N_i |k} 
R(\lambda,y, x_{N_1}, \ldots, x_{N_b}) =
Q(\lambda,y, \underline{x})
$$
where $Q(\lambda,y, \underline{x})$ is defined in \eqref{gform}.
Hence i) is proved. 

As for \eqref{fine3}, it suffices to show that the following identity
of generating series holds:
\begin{equation}
\label{realfine}
 - \sum_{g \geq 2} e({\overline{\mathcal M}_g}) \log\left(1 -
\lambda^{2g-2}\right) = \sum_{g \geq 2} u_g \lambda^{2g-2}.
\end{equation} 

$e({\overline{\mathcal M}_g})$ is the sum of the Euler characteristics
$e(\Delta(G))$, where $\Delta(G)$ are the strata in
${\overline{\mathcal M}}_g$. Analogously to Section \ref{sec6},
$e(\Delta(G))$ can be computed by taking into account connected
coverings of Riemann surfaces of genus $g$. 

Let us expand, now, the left hand side of \eqref{realfine}. Clearly, we get
\begin{equation}
\label{verafine} 
\sum_{g \geq 2} e({\overline{\mathcal M}}_g)\sum_{m \geq 1}
\frac{{\left(\lambda^{(2g-2)}\right)^m}}{m}.
\end{equation}

From what recalled on $e({\overline{\mathcal M}_g})$, \eqref{verafine}
can be interpreted as a generating series for coverings with more than
one connected component. This is exactly what the numbers $u_g$ enumerate.

\smallskip

\centerline{\sc Table 2: Some values of $e({\overline{\mathcal M}}_g^n)$}

\begin{center}
\begin{tabular}{r|r|r|r|r|r|r|r|}
$g \Big\backslash n$ & $0$ & $1$ & $2$ & $3$ & $4$ & $5$ & $6$ \\
 \hline 
$2$ & $6$ & $13$ & $42$ & $181$ & $1004$ & $6883$ &
 $56392$ \\ \hline
$3$ & $32$ & $102$ & $454$ & $2612$ & $18515$ & $156094$ &
 $1526677$ \\ \hline $4$ & $200$ & $882$ & $5214$ & $37945$ & $327584$
 & $3272624$ & $37151502$ \\ \hline
\end{tabular}
\end{center}

For genus $g=0,1$ our numbers coincide with the known values. For
$g=2$ our method showed an incogruence with the values in \cite{bgm}. In what follows, we adopt the same notation as in that paper. By Proposition 3.15, p. 507 in \cite{bgm}, the contribution of {\em Graphs of type 5} should be 
$$
\frac{1}{24(1-E)(1+D)^2}+\frac{11+2D-3D^2}{24(1-E)}+ \frac{1}{2}+ \frac{3D}{2} + \frac{7D^2}{4} + \frac{7D^3}{6} + \frac{11D^4}{24} + \frac{D^5}{8} + \frac{D^6}{48}
$$
and not
$$
\frac{1}{24(1-E)(1+D)^2}+\frac{11+2D-3D^2}{24(1-E)}+ \frac{1}{2}+ \frac{3D}{2} + \frac{7D^2}{4} + \frac{7D^3}{6} + \frac{11D^4}{24} - \frac{D^5}{8} + \frac{D^6}{48}.
$$

As a consequence, the final generating function $K_2(t)$ in \cite{bgm}
should be modified by adding $D^5/4$. This yields the same values we
get in the present paper for $g=2$.

\bigskip

{\bf Acknowledgements.}  Much of this work was performed while the
first author was a student and the second a visitor at the Scuola
Normale Superiore in Pisa.  Both authors would like to thank Professor
Enrico Arbarello, not only for organizing the special year 97-98 at
the institute, but also for his personal support and undying passion
for the moduli spaces.  He has inspired so many of us.

The first author would also like to thank Duke University for hosting
his visit in Fall 1998.

\end{document}